\documentclass[11pt]{amsart}
\usepackage{amsfonts,amssymb,amscd,amsmath,enumerate,verbatim,calc,amsthm}
\input xy
\xyoption{all}

\headsep=1cm \numberwithin{equation}{section}
\newtheorem{thm}{Theorem}[section]
\newtheorem{cor}[thm]{Corollary}
\newtheorem{lem}[thm]{Lemma}
\newtheorem{prop}[thm]{Proposition}
\newtheorem{defn}[thm]{Definition}
\newtheorem{hy}[thm]{Hyperhomology}

\newtheorem{exam}[thm]{Example}
\newtheorem{rem}[thm]{Remark}

\newtheorem*{theorem*}{Theorem}

\newcommand{\Hom}{\mbox{Hom}\,}

\newcommand{\Ext}{\mbox{Ext}\,}
\newcommand{\Tor}{\mbox{Tor}\,}
\newcommand{\Spec}{\mbox{Spec}\,}

\newcommand{\Supp}{\mbox{Supp}\,}

\newcommand{\gr}{\mbox{grade}\,}

\newcommand{\vdim}{\mbox{vdim}\,}
\newcommand{\depth}{\mbox{depth}\,}
\renewcommand{\dim}{\mbox{dim}\,}

\newcommand{\cd}{\mbox{cd}\,}

\newcommand{\T}{\mbox{T}\,}

\newcommand{\pd}{\mbox{pd}\,}
\newcommand{\id}{\mbox{id}\,}
\newcommand{\fd}{\mbox{fd}\,}
\newcommand{\gd}{\mbox{G-dim}\,}

\newcommand{\h}{\mbox{ht}\,}

\newcommand{\E}{\mbox{E}}

\newcommand{\uhom}{{\mathbf R}\Hom}
\newcommand{\utp}{\otimes^{\mathbf L}}

\renewcommand{\H}{\mbox{H}}
\newcommand{\V}{\mbox{V}}

\newcommand{\D}{\mbox{D}}

\newcommand{\fa}{\mathfrak{a}}
\newcommand{\fb}{\mathfrak{b}}

\newcommand{\fm}{\mathfrak{m}}

\newcommand{\fc}{\mathfrak{c}}

\begin{document}
\title[Duality for local cohomology modules]
 {some duality and equivalence results }

 \author[M.R. Zargar]{Majid Rahro Zargar }

\address{Faculty of Mathematical Sciences Department of Mathematics
University of Mohaghegh Ardabili 56199-11367, Ardabil, Iran.}

\email{zargar9077@gmail.com}
\email{m.zargar@uma.ac.ir}



\maketitle
\begin{abstract}Let $(R,\fm)$ be a relative Cohen-Macaulay local ring with respect to an ideal $\fa$ of $R$ and set $c:=\h\fa$. In this paper, we investigate some properties of the Matlis dual $\H_{\fa}^c(R)^{\vee}$ of the $R$-module $\H_{\fa}^c(R)$ and we show that such modules treat like canonical modules over Cohen-Macaulay local rings. Also, we provide some duality and equivalence results with respect to the module $\H_{\fa}^c(R)^{\vee}$ and so these results lead to achieve generalizations of some known results, such as the Local Duality Theorem, which have been provided over a Cohen-Macaulay local ring which admits a canonical module.
\end{abstract}
\maketitle
\tableofcontents
\section{Introduction}
Throughout this paper, $R$ is a commutative Noetherian ring, $\fa$ is a proper ideal of $R$ and $M$ is an $R$-module. In the case where $R$ is local with maximal ideal $\fm$, $\hat{R}$ denotes the $\fm$-adic completion of $R$, $\E(R/\fm)$ denotes the injective hull of the residue field $R/\fm$ and $^{\vee}$ denotes the Matlis dual functor $\Hom_{R}(-,\E(R/\fm))$. The theory of canonical modules for Cohen-Macaulay local rings is developed by Bruns and Herzog in \cite[Chapter 3]{HERZ}. But, in the more general setting of an arbitrary $n$-dimensional local ring $(R,\fm)$, a canonical module for $R$ is a finitely generated $R$-module $C$ such that $C\otimes_{R}\widehat{R}\cong\H_{\fm}^n(R)^{\vee}$. In the special case where R is Cohen.Macaulay, this condition turns out to be equivalent to Bruns' and Herzog's definition, namely that a canonical module is precisely a maximal Cohen-Macaulay module of type one and of finite injective dimension that in some of the literature is called a dualizing module. As a remarkable result, Foxby \cite{FOX}, Reiten \cite{RET} and Sharp \cite{SHARP} proved that a Cohen-Macaulay local ring $R$ admits a canonical module if and only if it is homomorphic image of a Gorenstein local ring. In particular, if $R$ is complete and Cohen-Macaulay, then $\H_{\fm}^n(R)^{\vee}$ is a canonical module of $R$. On the other hand, the Local Duality Theorem provides a fundamental tool for the study of local cohomology modules with respect to the maximal ideal of a local ring. Although it only applies to local rings which can be expressed as homomorphic images of Gorenstein local rings, this is not a great restriction, because the class of such local rings includes the local rings of points on affine and quasi-affine varieties, and, as mentioned above, all complete local rings. This theorem provides a functorial isomorphism between the Ext modules of the canonical module and the local cohomology modules with respect to the maximal ideal of a local ring. Over a Cohen-Macaulay local ring, a canonical module, if exists, plays an important role in the studying of the algebraic and homological properties of ideals and modules. Thus, finding modules which preserve beneficiaries of the canonical modules is the aim of many commutative algebraists. In this direction,
the principal aim of this paper is to study the properties of the $R$-module $\D_{\fa}:=\H_{\fa}^c(R)^{\vee}$ in the case where $R$ is a relative Cohen-Macaulay local ring with respect to $\fa$ and $c=\h_{R}\fa$ (i.e. there is precisely one non-vanishing local cohomology module of $R$ with
respect to $\fa$) and to provide a connection between the module $\H_{\fa}^c(R)^{\vee}$ and the local cohomology module with respect to the ideal $\fa$ of $R$. Indeed, we show that these modules treat like canonical modules over Cohen-Macaulay local rings. Recently, such modules have been studied by some authors, such as Hellus and Str$\ddot{u}$ckrad \cite{HST}, Hellus \cite{H}, Hellus and Schenzel \cite{HSC}, Khashyarmanesh \cite{KHSH} and Schenzel \cite{SC}, and has led to some interesting results.

The organization of this paper is as follows. In section 2, we collect some notations and definitions which will be used in the present paper. In section 3, first as a generalization of the Local Duality Theorem, we provide the following result:
\begin{thm}Let $\fa$ be an ideal of $R$, $E$ be an injective $R$-module and let $Y$ be an arbitrary complex of $R$-modules such that \emph{$\H_{\fa}^i(Y)=0$} if and only if $i\neq c$ for some integer $c$. Then, for all integer $i$ and for all complexes of $R$-modules $X$, there are the following isomorphisms:
 \begin{itemize}
\item[{(i)}]\emph{$\Tor^{R}_{c-i}(X,\H^{c}_{\fa}(Y))\cong\H_{\fa}^{i}(X\utp_R Y)$.}
\item[(ii)]\emph{$\Ext_{R}^{c-i}(X,\Hom_{R}(\H^{c}_{\fa}(Y),E))\cong\Hom_{R}(\H_{\fa}^i(X\utp_R Y),E)$.}
\end{itemize}
\end{thm}
Also, among of other things as an application of the above theorem in Proposition 3.4, we provide a connection between the module $\D_{\fa}$ and the local cohomology modules with respect to the ideal $\fa$.
Next, in Theorem 4.3, we state one of the main results which shows that the module $\D_{\fa}$ has some of the basic properties of ordinary canonical modules over Cohen-Macaulay local rings. Indeed, we prove the following result:
 \begin{thm}Let $(R,\fm)$ be an \emph{$\fa$-RCM} local ring with \emph{$\h_{R}\fa=c$}. Then the following statements hold true.
  \begin{itemize}
\item[(i)]{\emph{For all ideals $\fb$ of $R$ such that $\fa\subseteq\fb$, $\H_{\fb}^i(\D_{\fa})=0$ if and only if $i\neq c$. }}
\item[(ii)]{\emph{$\id_{R}(\D_{\fa})=c.$}}
\item[(iii)]{\emph{$\mu^{c}(\fm,\D_{\fa})=1.$}}
\item[(iv)]{\emph{For all $t= 0, 1,\ldots,c $ and for all $\fa$-RCM $R$-modules $M$ with $\cd(\fa,M)=t$, one has
\begin{itemize}
\item[(a)]{{$\Ext_{R}^{i}(M,\D_{\fa})=0$ if and only if $i\neq c-t$.}}
\item[(b)]{{$\H_{\fa}^i(\Ext_{R}^{c-t}(M,\D_{\fa}))=0$ if and only if $i\neq t$}}
\end{itemize}
}}
\end{itemize}
\end{thm}

 On the other hand, in Lemma 4.5, it is shown that {$\H_{\fb}^c(\D_{\fa}(R))\cong\E(R/\fm)$} for all ideals $\fb$ such that $\fa\subseteq\fb$. Also, in Proposition 4.8, we provide a result which leads us to determine the endomorphism ring of $\D_{\fa}$ and its non-zero local cohomology module at support $\fa$. Next, in section 5, first we introduce the concept of maximal relative Cohen-Macaulay modules with respect to an ideal $\fa$ of $R$. Then, in Proposition 5.1, we generalize a result and, in Remark 5.3, we raise a quite long-standing known conjecture related to the maximal Cohen-Macaulay modules.
 In addition, in Theorem 5.4, we generalize another known result related to ordinary canonical modules of Cohen-Macaulay local rings. Indeed, we show that if $R$ is a relative Cohen-Macaulay local ring with respect to $\fa$, then for all {$\fa$-RCM} $R$-modules $M$ with $\cd(\fa,M)=t$, there exists a natural isomorphism $M\otimes_{R}\hat R\cong\Ext^{c-t}_{R}(\Ext^{c-t}_{R}(M,\D_{\fa}),\D_{\fa}).$ In particular, if $M$ is maximal $\fa$-RCM, then $M\otimes_{R}\hat R\cong\Hom_{R}(\Hom_{R}(M,\D_{\fa}),\D_{\fa})$. Also, in Theorem 5.5, over relative Cohen-Macaulay local rings, we establish a characterization of maximal relative Cohen-Macaulay modules.

Now, suppose for a moment that $R$ is a Cohen-Macaulay local ring with a dualizing module $\Omega_{R}$ and that $M$ is a non-zero finitely generated $R$-module. Kawasaki, in \cite[Theorem 3.1]{KAW}, showed that if $M$ has finite projective dimension, then $M$ is Cohen-Macaulay if and only if $M\otimes_{R}\Omega_{R}$ is Cohen-Macaulay. Next, in \cite[Theorem 1.11]{KHY}, Khatami and Yassemi generalized this result. Indeed, they showed that the above result hold true whenever $M$ has finite Gorenstein dimension. In Theorem 5.6, we establish another main result which, by using $\D_{\fa}$ instead of $\Omega_{R}$, provides a generalization of the above mentioned result of Khatami and Yassemi. Finally, in Proposition 5.9 as an application of the above theorem, we could find another generalization of the results \cite[Theorem 3.1]{KAW} and \cite[Theorem 1.11]{KHY}.

\section{Notation and Prerequisites}
\begin{hy}\emph{ An $R$-complex $X$ is a sequence of $R$-modules $(X_v)_{v\in\mathbb{Z}}$ together
with $R$-linear maps $(\partial_v^{X}:X_v \longrightarrow X_{v-1})_{v\in\mathbb{Z}}$,
$$X =\cdots\longrightarrow X_{v+1}\stackrel{\partial_{v+1}^X} \longrightarrow X_{v}\stackrel{\partial_{v}^X}\longrightarrow X_{v-1}\longrightarrow\cdots,$$
such that $\partial_{v}^{X}\partial_{v+1}^{X}=0$ for all $v\in\mathbb{Z}$. For any $R$-module $M$, $\Gamma_{\fa}(M)$ is defined by
{$\Gamma_{\fa}(M):=\{x\in M~|~\Supp_R(Rx)\subseteq \V(\fa)\}.$} The right derived
functor of the functor $\Gamma_{\fa}(-)$ exists in $\mathrm{D}(R)$ and the complex ${\bf
R}\Gamma_{\fa}(X)$ is defined by ${\bf R}
\Gamma_{\fa}(X):=\Gamma_{\fa}(I)$, where $I$ is any
injective resolution of the complex $X$.  Then, for any
integer $i$, the $i$-th local cohomology module of $X$ with respect
to $\fa$ is defined by $\H_{\fa}^i(X):=\H_{-i}({\bf R}
\Gamma_{\fa}(X))$. Also, there is the functorial isomorphism ${\bf{R}}\Gamma_{\fa}(X)\simeq\check{C}_{\b{x}}\otimes_R^{{\bf L}}{X}$ in the derived category $\mathrm{D}(R)$, where ${\check{C}}_{\b{x}}$ denotes the $\check{C}$ech complex of $R$ with respect to $\b{x}=x_{1},\ldots, x_{t}$ and $\sqrt{\fa} = \sqrt{(x_{1},\ldots,x_{t})}$. Hence, this implies that $\H_{\fa}^i(X)\cong\H_{-i}{(\check{C}_{\b{x}}\otimes_R^{{\bf L}}{X})}$(for more details refer to \cite{LIP}). Let ${X}\in \mathrm{D}(R)$ and/or ${Y}\in \mathrm{D}(R)$.
The left-derived tensor product complex of ${X}$ and ${Y}$ in $\mathrm{D}(R)$ is
denoted by ${X}\otimes_R^{{\bf L}}{Y}$ and is defined by
{$${X}\otimes_R^{{\bf L}}{Y}\simeq {F}\otimes_R{Y}\simeq{X}
\otimes_R{F}^{'}\simeq {F}\otimes_R{F},^{'}$$} where ${F}$ and
${F}^{'}$ are flat resolutions of ${X}$ and ${Y}$, respectively.
Also, let ${X}\in \mathrm{D}(R)$ and/or ${Y}\in \mathrm{D}(R)$. The right derived homomorphism complex of ${X}$ and ${Y}$ in
$\mathrm{D}(R)$ is denoted by ${\bf R}\Hom_R({X},{Y})$ and is defined by $${\bf R}\Hom_R({X},{Y})\simeq \Hom_R({P},{Y})\simeq
\Hom_R({X},{I})\simeq \Hom_R({P},I),$$ where ${P}$ and
${I}$ are projective resolution of ${X}$ and injective resolution of ${Y}$,
respectively. For any two complex $X$ and $Y$, we set $\Ext_R^i(X,Y)=\H_{-i}(\uhom_R(X,Y))$ and $\Tor_{i}^R(X,Y)=\H_{i}(X\utp_{R}Y)$. For any integer $n$, the $n$-fold shift of a complex $(X,\xi^X)$ is the complex $\Sigma^nX$ given by $(\Sigma^n X)_v=X_{v-n}$ and $\xi_{v}^{\Sigma^nX}=(-1)^n\xi_{v-n}^{X}$. Also, we have $\H_i(\Sigma^nX)=\H_{i-n}(X)$. Next, for any contravariant, additive and exact functor $\T: C(R)\longrightarrow C(R)$ and $X\in C(R)$, where $C(R)$ denotes the category of $R$-complexes, we have $\H_l(\T(X))\cong\T(\H_{-l}(X))$.}
\end{hy}
\begin{defn} \emph{We say that a finitely generated $R$-module $M$ is relative Cohen-Macaulay with respect to $\fa$ if there is precisely one non-vanishing local cohomology module of $M$ with respect to $\fa$. Clearly this is the case if and only if $\gr(\fa,M)=\cd(\fa,M)$, where $\cd(\fa,M)$ is the largest integer $i$ for which $\H_{\fa}^i(M)\neq0$. For the convenance, we use the notation $\fa$-RCM, for an $R$-module which is relative Cohen-Macaulay with respect to $\fa$. Furthermore, a non-zero finitely generated $R$-module $M$ is said to be maximal $\fa$-RCM if $M$ is $\fa$-RCM with {$\cd(\fa,M)=\cd(\fa,R).$}}
\end{defn}
Observe that the above definition provides a generalization of the concept of Cohen-Macaulay and maximal Cohen-Macaulay modules. Also, notice that the notion of relative Cohen-Macaulay modules is connected with the notion of cohomologically complete intersection ideals which has been studied in \cite{HSC1} and has led to some interesting results. Furthermore, recently, such modules have been studied in \cite{HSC}, \cite{MZ} and \cite{MZ1}.
\begin{defn}\emph{For any $R$-module $M$ and a proper ideal $\fa$ of $R$, we set $\gr(\fa,M):=\inf\{i |~\Ext_{R}^i(R/\fa,M)\neq0\}.$ Notice that $\gr(\fa,M)$ is the least integer $i$ such that $\H_{\fa}^i(M)\neq0$. For any local ring $(R,\fm)$, we define $\D_{\bullet}(\bullet):=\Hom_{R}(\H_{\bullet}^{\tiny{\gr(\bullet,\bullet)}}(\bullet),\E(R/\fm))$. In the case where $R$ is $\fa$-RCM, for the convenience we set $\D_{\fa}:=\D_{\fa}(R)$.}
\end{defn}
\section{Generalized local duality}

The starting point of this section is the following result which plays an essential role in the present paper. This result provides a generalization of the local duality theorem. Indeed the generalization is done for ideals $\fa$ of the given ring $R$ and complexes $Y$ which their local cohomology modules with respect to $\fa$ have precisely one non-vanishing. The generalization is such that local duality becomes the special case where $Y=R$ and $\fa$ is the maximal ideal $\fm$ of the given local ring \emph{$(R,\fm)$}.
\begin{thm}Let $\fa$ be an ideal of $R$, $E$ be an injective $R$-module and let $Y$ be an arbitrary complex of $R$-modules such that \emph{$\H_{\fa}^i(Y)=0$} if and only if $i\neq c$ for some integer $c$. Then, for all integer $i$ and for all complexes of $R$-modules $X$, there are the following isomorphisms:
 \begin{itemize}
\item[{(i)}]\emph{$\Tor^{R}_{c-i}(X,\H^{c}_{\fa}(Y))\cong\H_{\fa}^{i}(X\utp_R Y)$.}
\item[(ii)]\emph{$\Ext_{R}^{c-i}(X,\Hom_{R}(\H^{c}_{\fa}(Y),E))\cong\Hom_{R}(\H_{\fa}^i(X\utp_R Y),E)$.}
\end{itemize}
\begin{proof}(i). Let $\b{x}=x_{1},\ldots,x_{t}$ be elements of $R$ such that $\sqrt{\fa} = \sqrt{(x_{1},\ldots,x_{t})}$ and let ${\check{C}}_{\b{x}}$ denotes the $\check{C}$ech complex of $R$ with respect to $x_{1},\ldots, x_{t}$. Now notice that $\H_{\fa}^i(Y)\cong\H_{-i}({\check{C}}_{\b{x}}\utp_{R}Y)$ for all $i$. Thus, by our assumption we have $\H_j({\check{C}}_{\b{x}}\utp_{R}Y)=0$ for all $j\neq -c$. Hence, one can deduce that there exists an isomorphism ${\check{C}}_{\b{x}}\utp_{R}Y\simeq\H({\check{C}}_{\b{x}}\utp_{R}Y)$ in derived category $\mathrm{D}(R)$, and so $\Sigma^{-c}\H_{\fa}^c(Y)\simeq {\check{C}}_{\b{x}}\utp_{R}Y$. Therefore, we can use the following isomorphisms
\[\begin{array}{rl}
 \Tor^{R}_{c-i}(X,\H^{c}_{\fa}(Y))&\cong\H_{c-i}(X\utp_R\H_{\fa}^c(Y))\\
&\cong\H_{-i}(\Sigma^{-c} (X\utp_R\H_{\fa}^c(Y)))\\
&\cong\H_{-i}(X\utp_R \Sigma^{-c}\H_{\fa}^c(Y))\\
&\cong\H_{-i}(X\utp_R {\check{C}}_{\b{x}}\utp_R Y)\\
&\cong\H_{\fa}^i(X\utp_R Y),
\end{array}\]
where the third isomorphism follows from \cite[Lemma 2.4.9]{CF1}, to compte the proof.

(ii). Note that since $\E$ is an injective $R$-module, $\uhom_{R}(\H^{c}_{\fa}(Y),E)\simeq\Hom_{R}(\H^{c}_{\fa}(Y),E)$. It therefore follows, by the following isomorphisms
\[\begin{array}{rl}
 \Ext_{R}^{c-i}(X,\Hom_{R}(\H^{c}_{\fa}(Y),E))&\cong\H_{i-c}(\uhom_R(X,\Hom_{R}(\H^{c}_{\fa}(Y),E)))\\
&\cong\H_{i-c}(\uhom_R(X,\uhom_{R}(\H^{c}_{\fa}(Y),E)))\\
&\cong\H_{i-c}(\uhom_R(X\utp_R\H^{c}_{\fa}(Y),E))\\
&\cong\Hom_R(\H_{c-i}(X\utp_R\H^{c}_{\fa}(Y)),E)\\
&\cong\Hom_R(\Tor_{c-i}^R(X,\H^{c}_{\fa}(Y)),E)\\
&\cong\Hom_R(\H_{\fa}^i(X\utp_R Y),E),
\end{array}\]
as required. Here notice that third isomorphism follows from \cite[Theorem 4.4.2]{CF1}.
\end{proof}
\end{thm}
The following result is the special case of the above theorem for category of $R$-modules.
\begin{prop}Let $\fa$ be an ideal of $R$, $E$ be an injective $R$-module and let $N$ be an arbitrary $R$-module such that \emph{$\H_{\fa}^i(N)=0$} if and only if $i\neq c$ for some integer $c$. Then, for all integers $i$ and for all $R$-modules $M$ such that \emph{$\Tor_j^R(M,N)=0$} for all $j>0$, there are the following isomorphisms:
 \begin{itemize}
\item[{(i)}]\emph{$\Tor^{R}_{c-i}(M,\H^{c}_{\fa}(N))\cong\H_{\fa}^{i}(M\otimes_R N)$.}
\item[(ii)]\emph{$\Ext_{R}^{c-i}(M,\Hom_{R}(\H^{c}_{\fa}(N),E))\cong\Hom_{R}(\H_{\fa}^i(M\otimes_R N),E)$.}
\end{itemize}
\begin{proof}Notice that since $\Tor^R_{j}(M,N)=0$ for all $j>0$, there is the isomorphism $M\otimes_R N\simeq M\utp_R N$ in the derived category. Hence, it immediately follows from Theorem 3.1.
\end{proof}
\end{prop}
As an consequence of the above result, we have the following result.
\begin{cor}Let $(R,\fm)$ be $d$-dimensional Cohen-Macaulay local ring with a dualizing $R$-module $\Omega_R$. Then, for all $R$-modules $M$ with finite Gorenstein flat dimension $($See \cite[Definition 4.5]{CF}$)$ there are the following isomorphisms:
 \begin{itemize}
\item[{(i)}]\emph{$\Tor^{R}_{d-i}(M,\E(R/\fm))\cong\H_{\fm}^{i}(M\otimes_R\Omega_R)$.}
\item[(ii)]\emph{$\Ext_{R}^{d-i}(M,\widehat{R})\cong\H_{\fm}^i(M\otimes_R \Omega_R)^{\vee}$.}
\end{itemize}
\begin{proof}First of all notice that $\Omega_R$ is a maximal Cohen-Macaulay $R$-module such that $\H_{\fm}^d(\Omega_R)\cong\E(R/\fm)$ and also, in view of \cite[Proposition 10.4.17]{EJ}, $\Tor_R^i(M,\Omega_R)=0$ for all $i>0$. Hence, one can use the known fact that $\widehat{R}\cong\Hom_{R}(\E(R/\fm),\E(R/\fm))$ and Proposition 3.2 to complete the proof.
\end{proof}
\end{cor}
The following result, which is a generalization of local duality theorem, is a consequence of Proposition 3.2 and recovers the \cite[Theorem 3.1]{HSC}.
\begin{prop}Let $(R,\fm)$ be an $\fa$-RCM local ring with \emph{$\h\fa=c$}. Then, for all $R$-modules $M$ and for all flat $R$-modules $F$  such that \emph{$\H_{\fa}^c(F)\neq 0$} $($e.g. faithfully flat $R$-modules$)$, there are the following isomorphisms:
 \begin{itemize}
\item[{(i)}]\emph{$\Tor^{R}_{c-i}(M,\H^{c}_{\fa}(F))\cong\H_{\fa}^{i}(F\otimes_R M)$.}
\item[(ii)]\emph{$\Ext_{R}^{c-i}(M,\D_{\fa}(F))\cong\H_{\fa}^i(F\otimes_R M)^{\vee}$.}
\end{itemize}
\begin{proof}First notice that, by \cite[Theorem 5.40]{ROT}, there is a directed index set $I$ and a family of finitely generated free $R$-modules $\{ F_{i}\}_{i\in I}$ such that $F=\underset{{i\in I}}\varinjlim F_{i}$. Notice that each $F_{i}$ is relative Cohen-Macaulay with respect to $\fa$ and that $\cd(\fa,{F_{i}})=c$. Therefore, $\H_{\fa}^{j}(F)=\underset{i\in I}\varinjlim\H_{\fa}^{j}(F_{i})=0$ for all $j\neq c$. Hence, one can use the Proposition 3.2 to complete the proof.
\end{proof}
\end{prop}
As an application of the Theorem 3.1, we prove the known Local Duality Theorem which has already been prove in \cite[Theorem 11.2.8]{BSH}.
 \begin{cor}Let $(R,\fm)$ be a $d$-dimensional Cohen-macaulay local ring with a dualizing $R$-module $\Omega_R$. Then, for all $R$-modules $M$ and all integers $i$, one has\emph{$$\H_{\fm}^i(M)\cong\Ext_{R}^{d-i}(M,\Omega_R)^{\vee}.$$ }
\end{cor}
\begin{proof} First notice that for all $i$, $\H_{\fm}^i(M)$ is Artinian, and so there is the $\hat R$-isomorphism and also $R$-isomorphism $\H_{\fm}^i(M) \cong\H_{\fm}^i(M)\otimes_R\hat R$. Therefore, using the Flat Base Change Theorem, the Matlis Duality Theorem \cite[Theorem 10.2.12]{BSH}, and the fact that $\D_{\fm\hat R}\cong\Omega_{R}\otimes_R\hat R$, implies that following isomorphisms:\[\begin{array}{rl}
\H_{\fm}^i(M)&\cong\H_{\fm\hat R}^i(\hat{M})\\
&\cong\Hom_{\hat R}(\Hom_{\hat R}(\H_{\fm\hat R}^i(\hat{M}),\E_{\hat R}(\hat R/\fm\hat R)),\E_{\hat R}(\hat R/\fm\hat R))\\
&\cong\Hom_{\hat R}(\Ext_{\hat R}^{d-i}(\hat{M},\D_{\fm\hat R}),\E_{\hat R}(\hat R/\fm\hat R))\\
&\cong\Hom_{\hat R}(\Ext_{\hat R}^{d-i}(\hat{M},\Omega_{R}\otimes_{R}\hat R),\E_{\hat R}(\hat R/\fm\hat R))\\
&\cong\Hom_{\hat R}(\Ext_{R}^{d-i}(M,\Omega_{R})\otimes_{R}\hat R,\E_{\hat R}(\hat R/\fm\hat R))\\
&\cong\Hom_{R}(\Ext_{R}^{d-i}(M,\Omega_{R}),\Hom_{\hat R}(\hat R,\E_{\hat R}(\hat R/\fm\hat R))\\
&\cong\Hom_{R}(\Ext_{R}^{d-i}(M,\Omega_{R}),\E_{\hat R}(\hat R/\fm\hat R))\\
&\cong\Ext_{R}^{d-i}(M,\Omega_{R})^{\vee}.
\end{array}\]
Notice that the third isomorphism follows from Proposition 3.4 and the last isomorphism follows from \cite[10.2.10]{BSH}.
\end{proof}
Notice that for a local Cohen-Macaulay ring $R$ with a dualizing $R$-module $\Omega_R$ and a finitely generated $R$-module $M$, we have $\sup\{~i~|~\Ext_{R}^i(M,\Omega_R)\neq 0\}=\depth R-\depth M.$ So, as a generalization of this result, we provide the following result which is an immediate consequence of the Proposition 3.4.
\begin{cor}Let $(R,\fm)$ be an \emph{$\fa$-RCM} local ring. Then, for all $R$-modules $M$, one has\emph{$$\sup\{~i~|~\Ext_{R}^i(M,\D_{\fa}(R))\neq 0\}=\gr(\fa,R)-\gr(\fa,M).$$ }
\end{cor}

\section{Generalized dualizing modules}
{Notice that if $(R,\fm)$ is a local complete ring with respect to $\fm$-adic topology, then $\D_{\fm}$ coincides to the ordinary dualizing module $\Omega_{R}$ of $R$. So the main aim in this section is to determine some properties of the module $\D_{\fa}$ and to show that how much such modules have behavior similar to the dualizing modules. In this direction, we need to the following two lemmas will play essential role in the proof of some our results. }
\begin{lem}Let $M$ be an $R$-module and $n$ be a non-negative integer. Then the following conditions are equivalent:
 \begin{itemize}
\item[{(i)}]\emph{{$\H^{i}_{\fa}(M)=0$ for all $i<n$.}
\item[(ii)]{$\Ext_{R}^i(N,M)=0$ for all $i<n$ and for any $\fa$-torsion $R$-module $N.$}
\item[(iii)]{$\Ext_{R}^i(R/\fa,M)=0$ for all $i<n$}.}
\end{itemize}
\begin{proof} (i)$\Rightarrow$(ii). Let $N$ be an $\fa$-torsion $R$-module. Then, it is straightforward to see that $\Hom_{R}(N,M)=\Hom_{R}(N,\Gamma_{\fa}(M))$. Hence, in view of [16, Theorem 10.47], we obtain the Grothendieck third quadrant spectral sequence with $$\E^{p,q}_{2}=\Ext_{R}^{p}(N, \H^{q}_{\fa}(M))\underset{p}\Longrightarrow \Ext_{R}^{p+q}(N,M).$$ Note that $E^{i-q,q}_{2}=0$ for all $0\leq i<n$ and $0\leq q\leq i$. Therefore $E^{i-q,q}_{\infty}=0$ for all $0\leq i<n$ and $0\leq q\leq i$. Now, consider the following filtration $$\{0\}=\Psi_{i+1}H^{i}\subseteq \Psi_{i}H^{i}\subseteq \cdots\subseteq\Psi_{-1}H^i\subseteq\Psi_{0}H^i=H^i,$$ where $H^i=\Ext_{R}^i(N,M)$ and $\E^{p,q}_{\infty}=\frac{\Psi_{p}H^i}{\Psi_{p+1}H^i}$, to see that
$\Ext_{R}^i(N,M)=0$ for all $i<n$. The implication (ii)$\Rightarrow$(iii) is trivial and the implication (iii)$\Rightarrow$(i) follows by a similar argument as the above.
\end{proof}
\end{lem}
The following lemma, which is needed to prove some next results, was proved in \cite[Propositon 3.1]{MZ1} by using the spectral sequence tools. For the convenience, we provide a new proof by using the derived category tools.
\begin{lem}Let $\fa$ be an ideal of $R$, $M$ be an $R$-module such that \emph{$\H_{\fa}^i(M)=0$} if and only if $i\neq c$ for some non-negative integer $n$. Then, for all $\fa$-torsion $R$-modules $N$ and for any integer $i$, one has \emph{$\Tor^{R}_{i}(N,\H_{\fa}^c(M))\cong\Tor_{i-c}^R(N,M).$}
\begin{proof}Let $N$ be an $\fa$-torsion $R$-module and let $\b{x}=x_{1},\ldots,x_{t}$ be elements of $R$ such that $\sqrt{\fa} = \sqrt{(x_{1},\ldots,x_{t})}$ and let ${\check{C}}_{\b{x}}$ denotes the $\check{C}$ech complex of $R$ with respect to $x_{1},\ldots, x_{t}$. Now notice that $\H_{\fa}^i(M)\cong\H_{-i}({\check{C}}_{\b{x}}\utp_{R}M)$ for all $i$. Then, by our assumption we have $\H_i({\check{C}}_{\b{x}}\utp_{R}M)=0$ for all $i\neq -c$. Hence, one can deduce that there exists an isomorphism ${\check{C}}_{\b{x}}\utp_{R}M\simeq\H({\check{C}}_{\b{x}}\utp_{R}M)$ in derived category, and so $\Sigma^{-c}\H_{\fa}^c(M)\simeq {\check{C}}_{\b{x}}\utp_{R}M$. Also, by similar argument, one has $N\simeq{\check{C}}_{\b{x}}\utp_{R}N$. Therefore, we can use the following isomorphisms:
\[\begin{array}{rl}
 \Tor^{R}_{i}(N,\H^{c}_{\fa}(M))&\cong\H_{i}(N\utp_R\H_{\fa}^c(M))\\
&\cong\H_{i-c}(\Sigma^{-c} (N\utp_R\H_{\fa}^c(M)))\\
&\cong\H_{i-c}(N\utp_R \Sigma^{-c}\H_{\fa}^c(M))\\
&\cong\H_{i-c}(N\utp_R {\check{C}}_{\b{x}}\utp_R M)\\
&\cong\H_{i-c}^i({\check{C}}_{\b{x}}\utp_R N\utp_R M)\\
&\cong\H_{i-c}^i(N\utp_R M)\\
&\cong\Tor_{i-c}^R(N,M),
\end{array}\]
to complete the proof.
\end{proof}
\end{lem}

As a first main result in this section, we provide the following result which determines some basic properties of the module $\D_{\fa}$.
\begin{thm}Let $(R,\fm)$ be an \emph{$\fa$-RCM} local ring with \emph{$\h_{R}\fa=c$}. Then the following statements hold true.
  \begin{itemize}
\item[(i)]{\emph{For all ideals $\fb$ of $R$ such that $\fa\subseteq\fb$, $\H_{\fb}^i(\D_{\fa})=0$ if and only if $i\neq c$. }}
\item[(ii)]{\emph{$\id_{R}(\D_{\fa})=c.$}}
\item[(iii)]{\emph{$\mu^{c}(\fm,\D_{\fa})=1.$}}
\item[(iv)]{\emph{For all $t= 0, 1,\ldots,c $ and for all $\fa$-RCM $R$-modules $M$ with $\cd(\fa,M)=t$, one has
\begin{itemize}
\item[(a)]{{$\Ext_{R}^{i}(M,\D_{\fa})=0$ if and only if $i\neq c-t$.}}
\item[(b)]{{$\H_{\fa}^i(\Ext_{R}^{c-t}(M,\D_{\fa}))=0$ if and only if $i\neq t$}}
\end{itemize}
}}
\end{itemize}
\begin{proof}(ii). First notice that since $R$ is $\fa$-RCM, one can use \cite[Theorem 3.1]{MZ1} to see that $\fd_{R}(\H_{\fa}^c(R))=c$, and so $\id_{R}(\D_{\fa})=c$.

(i). Let $\fb$ be an ideal of $R$ such that $\fa\subseteq\fb$. Then, by part(ii), $\H_{\fb}^{i}(\D_{\fa})=0$ for all $i> c$. Now, we show that $\H_{\fb}^{i}(\D_{\fa})=0$ for all $i<c$. To this end, by Lemma 4.1, it is enough to show that $\Ext_{R}^{i}(R/\fb,\D_{\fa})=0$ for all $i<c$. But, by Proposition 3.4, we have $$\Ext_{R}^{i}(R/\fb,\D_{\fa})\cong\H_{\fa}^{c-i}(R/\fb)^{\vee}.$$
Hence the assertion is done because of $\fa$-torsionness of the module $R/\fb$.

(iii). Since $\mu^{c}(\fm,\D_{\fa})=\vdim\Ext^c_{R}(R/\fm,\D_{\fa})$, it easily follows from Proposition 3.4.

(iv)(a). It follows from Proposition 3.4 and the our assumption.

(iv)(b). First notice that $\Supp(\Ext_{R}^{c-t}(M,\D_{\fa}))\subseteq\Supp(M)$ and that there exists a directed index set $I$ and a family of finitely generated  submodules $\{N_{i}\}_{i\in I}$ of $\Ext_{R}^{c-t}(M,\D_{\fa})$ such that $\Ext_{R}^{c-t}(M,\D_{\fa})\cong\underset{{i\in I}}\varinjlim N_{i}.$ Now, since $\Supp(N_{i})\subseteq\Supp(M)$ for all $i\in I$, by \cite[Theorem 2.2]{DNT}, $\cd(\fa,N_{i})\leq t$ for all $i$. Therefore, $\H_{\fa}^i(\Ext_{R}^{c-t}(M,\D_{\fa}))=0$ for all $i>t.$ Next, we show that $\H_{\fa}^i(\Ext_{R}^{c-t}(M,\D_{\fa}))=0$ for all $i<t.$ To this end, by using Lemma 4.1, it is enough to show that $$\Ext_{R}^i(R/\fa,\Ext_{R}^{c-t}(M,\D_{\fa}))=0$$ for all $i<t.$ Now, since $M$ is $\fa$-RCM, in view of Lemma 4.2 we have $\Tor_{i-t}^R(R/\fa, M)\cong\Tor_{i}^R(R/\fa,\H_{\fa}^{t}(M))$ for all $i$. Therefore, one can use Proposition 3.4 to get the following isomorphisms
\[\begin{array}{rl}
\Ext_{R}^i(R/\fa,\Ext_{R}^{c-t}(M,\D_{\fa}))&\cong\Tor_{i}^R(R/\fa,\H_{\fa}^{t}(M))^{\vee}\\
&\cong\Tor_{i-t}^R(R/\fa,M)^{\vee},
\end{array}\]
which complete the proof.
\end{proof}
\end{thm}
The following corollary which is some known results about a dualizing module is a consequence of the previous theorem.
\begin{cor}Let $(R,\fm)$ be $d$-dimensional Cohen-Macaulay local ring with a dualizing $R$-module $\Omega_R$. Then the following statements hold true.
  \begin{itemize}
\item[(i)]{\emph{$\H_{\fm}^i(\Omega_R)=0$ if and only if $i\neq d$. }}
\item[(ii)]{\emph{$\id_{R}(\Omega_R)=d.$}}
\item[(iii)]{\emph{$\mu^{d}(\fm,\Omega_R)=1.$}}
\item[(iv)]{\emph{For all $t= 0, 1,\ldots,d $ and for all Cohen-Macaulay $R$-modules $M$ with $\dim M=t$, one has
\begin{itemize}
\item[(a)]{{$\Ext_{R}^{i}(M,\Omega_R)=0$ if and only if $i\neq d-t$.}}
\item[(b)]{{$\H_{\fm}^i(\Ext_{R}^{d-t}(M,\Omega_R))=0$ if and only if $i\neq t$}}
\end{itemize}
}}
\end{itemize}
\begin{proof}It is straightforward to see that one may assume that $R$ is a complete local ring. Therefore, one has $\Omega_R\cong\D_{\fm}$; and so by using the Theorem 4.3 there nothing to prove.
\end{proof}
\end{cor}
The following lemma, which is needed in the proof of the next proposition, recovers \cite[Corollary 2.6]{KHSH}.
\begin{lem}Let the situation be as in Theorem 4.3. Then, for all ideals $\fb$ of $R$ such that $\fa\subseteq\fb$, there exists the isomorphism \emph{$\H_{\fb}^c(\D_{\fa})\cong\E(R/\fm)$.}
\end{lem}
\begin{proof} First, one can use Theorem 4.3(i)-(ii) and \cite[Theorem 2.5(i)]{MZ} to see that $\H_{\fb}^c(\D_{\fa})$  is injective. Now, in view of \cite[Proposition 2.1]{MZ} and Theorem 4.3(i), we get the following isomorphisms\[\begin{array}{rl}
\Hom_{R}(R/\fb,\H_{\fb}^c(\D_{\fa}))&\cong\Ext_{R}^c(R/\fb, \D_{\fa})\\
&\cong(R/\fb)^{\vee},
\end{array}\]
where the last isomorphism follows from Proposition 3.4(ii). Therefore, by Melkerson Lemma \cite[Theorem 7.2.1]{BSH}, $\H_{\fb}^c(\D_{\fa})$ is an Artinian $R$-module; and hence there exists a non-negative integer $t$ such that $\H_{\fb}^c(\D_{\fa})\cong\E(R/\fm)^t$. On the other hand, by \cite[Corollary 2.2]{MZ}, we have $\mu^0(\fm,\H_{\fb}^c(\D_{\fa}))=\mu^c(\fm, \D_{\fa})$. It therefore follows form Theorem 4.3(iii) that $t=1$.
\end{proof}
As an immediate consequence of the previous lemma we have the next result.
\begin{cor}Let $(R,\fm)$ be a $d$-dimensional local ring with \emph{$t:=\depth R$}. Let $x_1,\ldots,x_d$ be a system of parameters for $R$ which $x_1,\ldots,x_t$ is an $R$-sequence. Then, for all $0\leq i\leq t$ and for all $0\leq s\leq d-i$, we have \emph{$\H_{(x_1,\ldots,x_i)}^i(\D_{(x_1,\ldots,x_i)})\cong\H_{(x_1,\ldots,x_{i+s})}^i(\D_{(x_1,\ldots,x_i)})\cong\E(R/\fm)$.}
\end{cor}
Next, we provide an example to justify Theorem 4.3 and Lemma 4.5.
\begin{exam}\emph{Let $R=k[[x,y,z]]$ be the ring of power series, where $k$ is a field. Set $\fa=(x)$. Notice that $R$ is an $\fa$-RCM local ring with the maximal ideal $\fm=(x,y,z)$ and $\h_{R}{\fa}=1$. Consider the following exact sequence $$0\longrightarrow R\longrightarrow R_{x}\longrightarrow\H^1_{\fa}(R)\longrightarrow 0.$$ Now, by applying the functor $^{\vee}$ on the above exact sequence, one gets the exact sequence
\begin{equation}0\longrightarrow\D_{\fa}\longrightarrow R_{x}^{\vee}\longrightarrow\E(R/\fm)\longrightarrow 0. \end{equation}  Notice that $\Hom_{R}(R_{x},\E(R/\fm))$ is an injective $R$-module; and hence $\id_{R}(\D_{\fa})\leq1$. Now, we show that $\id_{R}(\D_{\fa})=1$. To this end, it is enough to show that $\D_{\fa}$ is not injective. Assume the contrary that it is injective. Since $\D_{\fa}\neq0$, there exists a non-zero element $f\in\D_{\fa}$ and a non-zero element $y\in\H_{\fa}^1(R)$ such that $0\neq f(y)\in\E(R/\fm)$; and hence there exists an element $r$ in $R$ such that $0\neq rf(y)\in R/\fm$. On the other hand, there exists an element $t\in\mathbb{N}$ such that $x^t y=0$. Now, since $\D_{\fa}$ is divisible, there exists an element $h\in\D_{\fa}$ such that $x^t h=f$, which is a contradiction. Hence $\id_{R}(\D_{\fa})=1$. Also, by (4.1), $\mu^{1}(\fm,\D_{\fa})=1$. Next, applying the functor $\Gamma_{\fa}(-)$ on (4.1) we obtain the exact sequence
\begin{equation}0\rightarrow\Gamma_{\fa}(\D_{\fa})\rightarrow \Gamma_{\fa}(R_{x}^{\vee})\rightarrow\Gamma_{\fa}(\E(R/\fm))\rightarrow\H^1_{\fa}(\D_{\fa})\rightarrow 0. \end{equation}
Therefore, since $\Gamma_{\fa}(R_{x}^{\vee})=0$, we have $\H_{\fa}^i(\D_{\fa})=0$ for all $i\neq1$ and $\H_{\fa}^1(\D_{\fa})\cong\E(R/\fm)$. Let $M=R/(y,z)$. Then, $\H_{\fa}^i(M)=0$ for all $i\neq1$. Since, $\id_{R}(\D_{\fa})=1$, $\Ext^i_{R}(M,\D_{\fa})=0$ for all $i>1$, Also, one can use (4.1) together with the fact that the natural map $M\longrightarrow M_{x}$ is one to one, to see that $\Ext^1_{R}(M,\D_{\fa})=0$. But, $\Hom_{R}(M,\D_{\fa})\cong\H_{\fa}^1(M)^{\vee}\neq0.$ Now, since $\Hom_{R}(R/\fa,\Hom_{R}(M,\D_{\fa}))\cong\H_{\fa}^1(M/\fa M)^{\vee}=0$, one has $\Gamma_{\fa}(\Hom_{R}(M,\D_{\fa}))=0$. On the other hand, one has $\H_{\fa}^i(\Hom_{R}(M,\D_{\fa}))=0$ for all $i>1$. Now, we must prove that $\H_{\fa}^1(\Hom_{R}(M,\D_{\fa}))\neq0$. To this end, by Lemma 4.1, it is enough to show that $\Ext_{R}^1(R/\fa,\Hom_{R}(M,\D_{\fa}))\neq0$. Since $$\Ext_{R}^1(R/\fa,\Hom_{R}(M,\D_{\fa}))\cong\Tor^R_{1}(R/\fa,\H_{\fa}^1(M))^{\vee},$$ we need only to show that $\Tor^R_{1}(R/\fa,\H_{\fa}^1(M))\neq0$. To this end, consider the exact sequence $0\longrightarrow M\longrightarrow M_{x}\longrightarrow\H^1_{\fa}(M)\longrightarrow 0$ to induce the exact sequence $$\Tor_{1}^R(R/\fa,\H_{\fa}^1(M))\rightarrow M/\fa M\xrightarrow{\varphi} M_{x}/\fa M_{x}\rightarrow\H^1_{\fa}(M)\otimes_{R}R/\fa\rightarrow 0,$$
and use the fact that $\varphi=0$ to complete the proof.}
\end{exam}

\begin{prop}Let $(R,\fm, k)$ be an \emph{$\fa$-RCM} local ring with \emph{$\h_{R}\fa=c$} and let \emph{$\fb\subseteq\fa$} be an ideal of $R$. Set \emph{$t:=\inf\{~ i~ | ~\H^{i}_\fb(\D_{\fa})\neq 0\}$} and consider the following statements:
\begin{itemize}
\item[(i)]{\emph{$\H^{i}_\fb(\D_{\fa})=0$ for all $i\neq t$.}}
\item[(ii)]{\emph{$\H^{c-t}_{\fa}(\H^{t}_{\fb}(\D_{\fa}))\cong \E_{R}(k)$ and $\H^{i}_{\fa}(\H^{t}_{\fb}(\D_{\fa}))=0$ for $i\neq c-t$ }.}
\item[(iii)]{\emph{ For $\fc=\fm$ and $\fc=\fa$, $\Ext_{R}^{c-t}(R/\fc,\H^{t}_\fb(\D_{\fa}))\cong\E_{R/\fc}(k)$ and $\Ext_{R}^{i}(R/\fc,\H^{t}_{\fb}(\D_{\fa}))=0$ for all $i\neq c-t$ }.}
\end{itemize}
Then, the implications \emph{(i)$\Rightarrow$(ii) and (ii)$\Leftrightarrow$(iii)} hold true. Moreover, if condition \emph{(i)} is satisfied, then \emph{$\Ext^{i}_{R}(\H^{t}_\fb(\D_{\fa}),\H^{t}_\fb(\D_{\fa}))=0$} for all $i\neq0$ and \emph{$\widehat{R}\cong\Hom_{R}(\H^{t}_\fb(\D_{\fa}),\H^{t}_\fb(\D_{\fa}))$.}
\end{prop}
\begin{proof}First, notice that one can use Lemma 4.1 and Proposition 3.4 to deduce that $t=c-\cd(\fa,R/\fb)$; and hence it is a finite number.

(i)$\Rightarrow$(ii). In view of \cite[Proposition 2.8]{MZ} and the assumption, one has the isomorphism  $\H^{i}_{\fa}(\H^{t}_{\fb}(\D_{\fa}))\cong\H^{i+t}_{\fa}(\D_{\fa})$ for all $i\geq0$. Therefore, one can use Theorem 4.3(i) and Lemma 4.5 to complete the proof.

(ii)$\Rightarrow$(iii): Let $\fc\in\{\fa,\fm\}$. Since $R/\fc$ is an $\fa$-torsion $R$-module, Lemma 4.1(i)$\Rightarrow$(ii) implies that $\Ext_{R}^{i}(R/\fc,\H^{t}_{\fb}(\D_{\fa}))=0$ for all $i<c-t$. On the other hand, since $\H^{c-t}_{\fa}(\H^{t}_{\fb}(\D_{\fa}))$ is injective, one can use \cite[Proposition 2.1]{MZ} to see that $\Ext_{R}^{i}(R/\fc,\H^{t}_{\fb}(\D_{\fa})))=0$ for all $c-t<i$ and that $\Ext_{R}^{c-t}(R/\fc,\H^{t}_{\fb}(\D_{\fa}))\cong\Hom_{R}(R/\fc,\H^{c-t}_{\fa}(\H^{t}_{\fb}(\D_{\fa})))$. Hence the proof is complete.

(iii)$\Rightarrow$(ii): First, in view of Lemma 4.1(iii)$\Rightarrow$(i) and our assumption for $\fc=\fa$, we deduce that $\H^{i}_{\fa}(\H^{t}_{\fb}(\D_{\fa}))=0$ for all $i<c-t$. Notice that, in view of Proposition 3.4, we have $$\Ext_{R}^{c-i}(\H_{\fb}^t(\D_{\fa}),\D_{\fa})\cong\H^{i}_{\fa}(\H^{t}_{\fb}(\D_{\fa}))^{\vee}$$ for all $i$. Now, since $\H^{i}_\fb(\D_{\fa})=0$ for all $i<t$ and $\H^{i}_\fb(\D_{\fa})$ is $\fb$-torsion, one can use Lemma 4.1(i)$\Rightarrow$(ii) to see that $\Ext_{R}^{j}(\H_{\fb}^t(\D_{\fa}),\D_{\fa})=0$ for all $j<t$; and hence $\H^{i}_{\fa}(\H^{t}_{\fb}(\D_{\fa}))=0 $ for all $i>c-t$.
Therefore, one can use \cite[Proposition 2.1]{MZ} and assumption for $\fc=\fa$ to deduce that \[\begin{array}{rl} \Hom_{R}(R/\fa,\H^{c-t}_{\fa}(\H^{t}_{\fb}(\D_{\fa})))&\cong\Ext_{R}^{c-t}(R/\fa,\H^{t}_\fb(\D_{\fa}))\\
&\cong(R/\fa)^{\vee}.
\end{array}\] Hence, $\Hom_{R}(R/\fa,\H^{c-t}_{\fa}(\H^{t}_{\fb}(\D_{\fa})))$ is an Artinian $R$-module. Therefore, by \cite[Theorem 7.1.2]{BSH}, $\H^{c-t}_{\fa}(\H^{t}_{\fb}(\D_{\fa}))$ is Artinian. Thus, one can use \cite[Corollary 2.2]{MZ} to see that $\mu^1(\fm,\H^{c-t}_{\fa}(\H^{t}_{\fb}(\D_{\fa})))=\mu^{c-t+1}(\fm,\H^{t}_{\fb}(\D_{\fa}))=0$. Hence $\H^{c-t}_{\fa}(\H^{t}_{\fb}(\D_{\fa}))$ is an injective $R$-module. Therefore, again we can use \cite[Corollary 2.2]{MZ} to see that $\H^{c-t}_{\fa}(\H^{t}_{\fb}(\D_{\fa}))\cong \E_{R}(k)$.

For the finial assertion, suppose that the statement (i) hold true. Then, we can use \cite[Proposition 2.1]{MZ} and Proposition 3.4 to get the following isomorphisms\[\begin{array}{rl}\Ext^{i}_{R}(\H^{t}_\fb(\D_{\fa}),\H^{t}_\fb(\D_{\fa}))&\cong\Ext_{R}^{i+t}(\H^{t}_\fb(\D_{\fa}),\D_{\fa})\\
&\cong\H^{c-i-t}_{\fa}(\H^{t}_{\fb}(\D_{\fa}))^{\vee}.
\end{array}\]
Now, one can use (ii) to complete the proof.
\end{proof}

The following corollary is an immediate consequence of the above proposition and Theorem 4.3. Notice that the first part of the next corollary has been proved in \cite[Theorem 1.3(b)]{SC}.
\begin{cor}Let $(R,\fm)$ be an \emph{$\fa$-RCM} local ring with \emph{$\h_{R}\fa=c$.} Then the following statements hold true: \begin{itemize}
\item[(i)]{\emph{$\hat{R}\cong\Hom_{R}(\H^{c}_\fa(\D_{\fa}),\H^{c}_\fa(\D_{\fa}))\cong\Hom_{R}(\D_{\fa},\D_{\fa})$.}}
\item[(ii)]{\emph{$\Ext_{R}^i(\D_{\fa},\D_{\fa})=\Ext^{i}_{R}(\H^{c}_\fa(\D_{\fa}),\H^{c}_\fa(\D_{\fa}))=0$ for all $i>0$.}}
\end{itemize}
\end{cor}

\section{Maximal relative Cohen-Macaulay modules}
Given an $R$-module $M$ and a finite free resolution \emph{$F=\{F_{i},\lambda_{i}^F\}$} of $M$, we define $\Omega_{i}^F(M)=\ker\lambda_{i-1}^F.$  The starting point of this section is the following proposition which recovers the well known fact that, for a finitely generated module $M$ over a Cohen-Macaulay local ring, $\Omega_{n}^F(M)$ is zero or maximal Cohen-Macaulay for all $n\geq\dim R$.
\begin{prop}Let $M$ be a finitely generated $R$-module and let $F$ be a finite free resolution for $M$. Then the following hold true:
  \begin{itemize}
\item[(i)]{\emph{$\gr(\fa,\Omega_{n}^F(M))\geq\min\{ n, \gr(\fa,R)\}$.}}
\item[(ii)]{\emph{If $R$ is $\fa$-RCM, then $\Omega_{n}^F(M)$ is zero or maximal $\fa$-RCM for all $n\geq\cd(\fa,R).$}}
\end{itemize}
\begin{proof}(i). We proceed by induction on $n$. If $n=0$, then there is nothing to prove. Let $n>0$ and suppose that the result has been proved for $n-1$. Let $$F= \cdots\longrightarrow F_{n}\longrightarrow F_{n-1}\longrightarrow\cdots\longrightarrow F_1\longrightarrow F_0\longrightarrow M\longrightarrow 0$$ be a finite free resolution for $M$. Then, there are the following exact sequences
$$0\longrightarrow\Omega_{n}^F (M)\longrightarrow F_{n-1}\longrightarrow\Omega_{n-1}^F(M)\longrightarrow 0$$
$$0\longrightarrow\Omega_{n-1}^F(M)\longrightarrow F_{n-2}\longrightarrow\cdots\longrightarrow F_1\longrightarrow F_0\longrightarrow M\longrightarrow 0.$$
Hence, one can use \cite[Proposition 1.2.9]{BH} and inductive hypothesis to see that { \[\begin{array}{rl}
\gr(\fa,\Omega_{n}^F(M))&\geq\min\{\gr(\fa, F_{n-1}), \gr(\fa,\Omega_{n-1}^F(M)) +1\}\\
&\geq\min\{\gr(\fa, F_{n-1}), \min\{ n-1, \gr(\fa,R)\} +1\}\\
&\geq\min\{\gr(\fa,R), n\}.
\end{array}\]}

(ii). Let $R$ be $\fa$-RCM and $n\geq\cd(\fa,R)$. Then, $\cd(\fa,R)=\gr(\fa,R).$ Next, since $\Supp(\Omega_{n}^F(M))\subseteq\Spec(R)$, by \cite[Theorem 2.2]{DNT}, we have $\cd(\fa,\Omega_{n}^F(M))\leq\cd(\fa,R).$ Therefore, one can use (i) to complete the proof.
\end{proof}
\end{prop}

An immediate application of the previous proposition is the next corollary.
\begin{cor}Let $(R,\fm)$ be a local ring, $x_{1},\ldots,x_{n}$ be an $R$-sequence and let $M$ be a non-zero finitely generated $R$-module. Then, for all $i\geq n$,  $x_{1},\ldots,x_{n}$ is an $\Omega_{i}^F M$-sequence, whenever $\Omega_{i}^F M$ is non-zero.
\end{cor}
\begin{rem} \emph{It is an important conjecture which indicates that any local ring admits a maximal Cohen-Macaulay module. We can raise the above conjecture for maximal relative Cohen-Macaulay modules. Indeed, we have the following conjecture.}
\end{rem}
\hspace{-0.45cm{}${\mathbf{Conjecture}}$. Let $\fa$ be a proper ideal of a local ring $R$. Then there exists a maximal relative Cohen-Macaulay module with respect to $\fa$ over $R$.}\\

Let $R$ be a Cohen-Macaulay local ring with a dualizing $R$-module $\Omega_{R}$. Then, there exists the known isomorphism {$M\cong\Ext^{d-t}_{R}(\Ext^{d-t}_{R}(M,\Omega_{R}),\Omega_{R})$}, whenever $M$ is a Cohen-Macaulay $R$-module of dimension $t$. In the following result, we use $\D_{\fa}$ instead of $\Omega_{R}$ to provide a generalization of the above result.

\begin{thm}Let $(R,\fm)$ be a local \emph{$\fa$-RCM} local ring with $\gr(\fa,R)=c$. Then, for all nonzero \emph{$\fa$-RCM} $R$-modules $M$ with $\cd(\fa,M)=t$, there exists a natural isomorphism\emph{$$M\otimes_{R}\widehat R\cong\Ext^{c-t}_{R}(\Ext^{c-t}_{R}(M,\D_{\fa}),\D_{\fa}).$$}In particular, if $M$ is maximal \emph{$\fa$-RCM}, then \emph{$M\otimes_{R}\widehat R\cong\Hom_{R}(\Hom_{R}(M,\D_{\fa}),\D_{\fa}).$}
\end{thm}
\begin{proof}  Let $M$ be an $\fa$-RCM $R$-module with $\cd(\fa,M)=t$. Then, in view of Theorem 4.3(iv), one has $\Ext_{R}^{j}(M,\D_{\fa})=0$ for all $j\neq c-t$. Therefore, for all $i$, one can use the following isomorphisms  \[\begin{array}{rl}
\Ext_{R}^i(M,\Sigma^{c-t}\D_{\fa})&\cong\H_{-i}(\uhom_{R}(M,\Sigma^{c-t}\D_{\fa}))\\
&\cong\H_{-i}(\Sigma^{c-t}\uhom_{R}(M,\D_{\fa}))\\
&\cong\H_{-i+t-c}(\uhom_{R}(M,\D_{\fa}))\\
&\cong\Ext_R^{i+c-t}(M,\D_{\fa}),
\end{array}\]
to deduce that $\Hom_{R}(M,\Sigma^{c-t}\D_{\fa})\cong\Ext_R^{c-t}(M,\D_{\fa})$ and $\Ext_{R}^i(M,\Sigma^{c-t}\D_{\fa})=0$ for all $i>0$. Notice that the second isomorphism follows from \cite[Lemma 2.3.10]{CF1}. Hence we obtain the isomorphism $\Hom_{R}(M,\Sigma^{c-t}\D_{\fa})\simeq \uhom_{R}(M,\Sigma^{c-t}\D_{\fa})$ in the derived category. On the other hand, in view of Corollary 4.9, we have $\uhom_R(\D_{\fa},\D_{\fa})\simeq\Hom_R(\D_{\fa},\D_{\fa})\cong\widehat{R}$. Hence, one can use Theorem 4.3(ii) and \cite[A.4.24]{CF0} to obtain the following isomorphism
\[\begin{array}{rl}\uhom_{R}(\uhom_{R}(M,\D_{\fa}),\D_{\fa})&\simeq M\utp_R\uhom_R(\D_{\fa},\D_{\fa})\\
&\simeq M\utp_R\widehat{R}.
\end{array}\]

Therefore, we can use the following isomorphisms
\[\begin{array}{rl}
\Ext_R^{c-t}(\Ext_{R}^{c-t}(M,\D_{\fa}),\D_{\fa})&\cong\H_{t-c}(\uhom_{R}(\Ext_{R}^{c-t}(M,\D_{\fa}),\D_{\fa}))\\
&\cong\H_{t-c}(\uhom_{R}(\Hom_{R}(M,\Sigma^{c-t}\D_{\fa}),\D_{\fa}))\\
&\cong\H_{t-c}(\uhom_{R}(\uhom_{R}(M,\Sigma^{c-t}\D_{\fa}),\D_{\fa}))\\
&\cong\H_{t-c}(\uhom_{R}(\Sigma^{c-t}\uhom_{R}(M,\D_{\fa}),\D_{\fa}))\\
&\cong\H_{t-c}(\Sigma^{t-c}\uhom_{R}(\uhom_{R}(M,\D_{\fa}),\D_{\fa}))\\
&\cong\H_{0}(\uhom_{R}(\uhom_{R}(M,\D_{\fa}),\D_{\fa}))\\
&\cong M\otimes_R\widehat{R},
\end{array}\]
to complete the proof. The forth isomorphism follows from \cite[Lemma 2.3.10]{CF1} and the fifth isomorphism follows from \cite[Lemma 2.3.16]{CF1}.
\end{proof}

The following theorem provides a characterization of maximal relative Cohen-Macaulay modules over a relative Cohen-Macaulay ring.
\begin{thm}Let $(R,\fm)$ be an \emph{$\fa$-RCM} local ring and let $M$ be a non-zero finitely generated $R$-module. Then the following statements are equivalent:
\begin{itemize}
\item[(i)]\emph{{$M$ is a maximal $\fa$-RCM} $R$-module.}
\item[(ii)]{\emph{$\Ext_{R}^i(M,\D_{\fa})=0$} for all $i>0$.}
\item[(iii)]\emph{The following conditions hold true:}
\begin{itemize}
\item[(a)]\emph{$M\otimes_{R}\hat R\cong\Hom_{R}(\Hom_{R}(M,\D_{\fa}),\D_{\fa}).$}
\item[(b)]{\emph{$\Ext_{R}^i(\Hom_{R}(M,\D_{\fa}),\D_{\fa})=0$} for all $i>0$.}
 \end{itemize}
\end{itemize}
\end{thm}
\begin{proof}The implications (i)$\Leftrightarrow$(ii) follows from Corollary 3.6 and the fact that $\cd(\fa,M)\leq\cd(\fa,R)$. Also, the implication (i)$\Rightarrow$(iii)(a) follows from Theorem 5.4.

(i)$\Rightarrow$(iii)(b). In view of Theorem 4.3(iv)(b) we have $\H_{\fa}^i(\Hom_{R}(M,\D_{\fa}))=0$ for all $i\neq\cd(\fa,M)$. Therefore, one can use Proposition 3.4 to complete the proof.

(iii)$\Rightarrow$(i). In view of Proposition 3.4, for all $i$, one has $$\Ext_{R}^i(\Hom_{R}(M,\D_{\fa}),\D_{\fa})\cong\H_{\fa}^{c-i}(\Hom_{R}(M,\D_{\fa}))^{\vee},$$ where $c=\cd(\fa,R)$. Therefore, by the assumption (iii)(b), one gets $\H_{\fa}^{j}(\Hom_{R}(M,\D_{\fa}))=0$ for all $j<c$. On the other hand, by the same argument as in the proof Theorem 4.3(iv)(b), one can see that $\cd(\fa,\Hom_{R}(M,\D_{\fa}))\leq\cd(\fa,R)$; and hence $\H_{\fa}^{j}(\Hom_{R}(M,\D_{\fa}))=0$ for all $j\neq c$. Therefore, one can use Proposition 3.4, the assumption (iii)(a) and Lemma 4.2 to obtain the following isomorphisms
\[\begin{array}{rl}
\Ext_{R}^i(R/\fa, M\otimes_{R}\hat R)&\cong\Ext_{R}^i(R/\fa,\H_{\fa}^c(\Hom_{R}(M,\D_{\fa}))^{\vee})\\
&\cong\Tor_{i}^R(R/\fa,\H_{\fa}^{c}(\Hom_{R}(M,\D_{\fa})))^{\vee}\\
&\cong\Tor_{i-c}^R(R/\fa,\Hom_{R}(M,\D_{\fa}))^{\vee}.
\end{array}\]
Hence, $\Ext_{R}^i(R/\fa, M\otimes_{R}\hat R)=0$ for all $i<c$, and so by Lemma 4.1, $\H_{\fa}^i(M\otimes_{R}\hat R)=0$ for all $i<c$. Therefore, using the Flat Base Change Theorem and Independence Theorem implies that $\H_{\fa}^i(M)=0$ for all $i<c$. Hence $c\leq\gr(\fa,M)$; and so $M$ is maximal $\fa$-RCM.
\end{proof}

 As we described in the introduction, the following theorem, which is one of the main results, provides a generalization of the result \cite[Theorem 1.11]{KHY} of Khatami and Yassemi. In the next theorem, we shall use the notion of grade of $M$ which is defined by $\gr(M)=\inf\{~i~|~\Ext_{R}^i(M,R)\neq0\}.$

\begin{thm}Let $(R,\fm)$ be an \emph{$\fa$-RCM} local ring of dimension $d$ with \emph{$\h_{R}\fa=c$} and let $M$ be a non-zero finitely generated $R$-module such that \emph{$\Tor^R_{i}(M,\D_{\fa})=0$} for all $i>0$. Then the following statements hold true: \begin{itemize}
\item[(i)]{\emph{$\Ext_{R}^i(M\otimes_{R}\D_{\fa},\D_{\fa})\cong\Ext^i_{R}(M,\hat{R})$ for all $i$ and $M\otimes_{R}\D_{\fa}\neq0$.}}
\item[(ii)]{Suppose that \emph{$\sup\{~ i~|~\Ext_{R}^i(M,R)\neq0\}=\depth R-\depth M$ and that $\gr(M)=\depth R-\dim M$. Then $M$ is Cohen-Macaulay if and only if $\H_{\fa}^{i}(M\otimes_{R}\D_{\fa})=0$ for all $i\neq c-\gr(M)$.}}
\end{itemize}
\end{thm}
\begin{proof}(i). Let $P$ be a projective $R$-module. Then, by \cite[Theorem 5.40]{ROT}, there exists a directed index set $I$ and a family of finitely generated free $R$-modules $\{ F_{j}\}_{j\in I}$ such that $P=\underset{{j\in I}}\varinjlim F_{j}$. Therefore, in view of Theorem 4.3(i), we have $\H_{\fa}^{c-i}(F_j\otimes_{R}\D_{\fa})=0$ for all $j\in I$ and for all $i\neq 0$, and so $\H_{\fa}^{c-i}(P\otimes_{R}\D_{\fa})\cong\underset{{i\in I}}\varinjlim\H_{\fa}^{c-i}(F_i\otimes_{R}\D_{\fa})=0$ for all $i\neq 0$. Hence, one can use Proposition 3.4 to see that $\Ext_{R}^i(P\otimes_{R}\D_{\fa},\D_{\fa})=0$ for all $i\neq0$. Hence, by using \cite[Theorem 10.49]{ROT} and Corollary 4.9(i), there is a third quadrant spectral sequence with $$\E_{2}^{p,q}:=\Ext_{R}^p(\Tor_{q}^R(M,\D_{\fa}),\D_{\fa})\underset{p}\Longrightarrow\Ext_{R}^{p+q}(M,\hat{R}).$$
Now, since {$\Tor^R_{i}(M,\D_{\fa})=0$} for all $i>0$, $E_{2}^{p,q}=0$ for all $q\neq0$. Therefore, this spectral sequence collapses in
the column $q=0$; and hence one gets, for all $i$, the isomorphism
$$\Ext_{R}^i(M\otimes_{R}\D_{\fa},\D_{\fa})\cong\Ext^i_{R}(M,\hat{R}),$$ as required. Next, notice that since $\widehat{R}$ is a flat $R$-module, by \cite[Theorem 3.2.15]{EJ}, $\Ext^i_{R}(M,\hat{R})\cong\Ext^i_{R}(M,R)\otimes_{R}\hat{R}$ for all $i$. Therefore, $M\otimes_{R}\D_{\fa}\neq0$.

(ii). Suppose that $M$ is Cohen-Macaulay. Then, $\Ext_{R}^i(M,R)=0$ for all $i\neq\gr(M)$. Hence, one can use (i) and Proposition 3.4 to complete the proof. The converse follows again by (i) and Proposition 3.4.
 \end{proof}

Next, we recall the concept of Gorenstein dimension which was introduced by Auslander in \cite{Aus}.
\begin{defn}\emph{A finite $R$-module $M$ is said to be of Gorenstein
dimension zero and we write {$\gd(M)=0$}, if and only if
\begin{itemize}
\item[(i)]{{$\Ext^i_{R}(M,R)=0$ for all $i>0$.}}
\item[(ii)]{{$\Ext^i_{R}(\Hom_{R}(M,R),R)=0$ for all $i>0$.}}
\item[(iii)]{{The canonical map $M\rightarrow\Hom_{R}(\Hom_{R}(M,R),R)$ is an isomorphism.}}
\end{itemize}
For a non-negative integer $n$, the $R$-module $M$ is said to be of
Gorenstein dimension at most $n$, if and only if there exists an exact
sequence
$$0\longrightarrow G_n\longrightarrow G_{n-1}\longrightarrow\cdots \longrightarrow G_0\longrightarrow M\longrightarrow0,$$
where {$\gd(G_{i})=0$} for $0\leq i\leq n$. If such a sequence does not exist, then we write {$\gd(N)=\infty$}.}
\end{defn}

Let $R$ be local and $M$ be a non-zero finitely generated $R$-module of finite $G$-dimension. Then \cite[Corollary 1.10]{CF} and
\cite[Proposition 10.4.17]{EJ} implies that {$\sup\{~ i~|~\Ext_{R}^i(M,R)\neq0\}=\depth R-\depth M$} and {$\Tor^R_{i}(M,\Omega_{R})=0$} for all $i>0$, where $\Omega_R$ is a dualizing $R$-module. Next, we provide a remark which shows that the converse of this fact is no longer true, that is, if {$\sup\{~ i~|~\Ext_{R}^i(M,R)\neq0\}=\depth R-\depth M$} and {$\Tor^R_{i}(M,\Omega_{R})=0$} for all $i>0$, then it is not necessary that $M$ has finite $G$-dimension.

\begin{rem}\emph{Let $(R,\fm)$ be a Cohen-Macaulay local ring which admits a dualizing module $\Omega_{R}$ and let $M$ be a finitely generated $R$-module such that {$\sup\{~ i~|~\Ext_{R}^i(M,R)\neq0\}=\depth R-\depth M$}. Next, suppose, in contrary, that the converse of the above mentioned fact is true. In this case, we prove that $M$ has finite $G$-dimension. To this end, we may assume that $\pd_{R}(M)=\infty$. Let $n\geq\dim R$ and let $\Omega_{n}^F(M)$ be the $n$-th syzygy of a finite free resolution $F$ for $M$. Then, in view of Proposition 5.1, $\Omega_{n}^F(M)$ is maximal Cohen-Macaulay. Also, since $\Ext_{R}^i(\Omega_{n}^F(M),R)\cong\Ext_{R}^{i+n}(M,R)$ for all $i>0$, one has $\Ext_{R}^i(\Omega_{n}^F(M),R)=0$ for all $i>0$. Therefore, $\Hom_{R}(\Omega_{n}^F(M),R)$ is maximal Cohen-Macaulay. Next, in view of \cite[A.4.24]{CF0}, we have the following isomorphisms:
\begin{equation} \textbf{R}\Hom_R(\textbf{R}\Hom_R(\Omega_{n}^F(M),R),\Omega_{R})\simeq \Omega_{n}^F(M)\otimes_{R}^{\textbf{L}}\textbf{R}\Hom_R(R,\Omega_{R})\simeq \Omega_{n}^F(M)\otimes_{R}^{\textbf{L}}\Omega_{R}
\end{equation}
in derived category. Since $\Ext_{R}^i(\Omega_{n}^F(M),R)=0$ for all $i>0$, we can deduce that $\textbf{R}\Hom_R(\Omega_{n}^F(M),R)\simeq\Hom_{R}(\Omega_{n}^F(M),R)$. On the other hand, in view of \cite[Corollary 3.5.11]{BH}, $\Ext_{R}^i(\Hom_R(\Omega_{n}^F(M),R),\Omega_{R})=0$ for all $i>0$. Therefore, by (5.1), $\Tor_{R}^{i}(\Omega_{n}^F(M),\Omega_{R})=0$ for all $i>0$. Thus, by using the contrary assumption and \cite[Lemma 1.2.6]{CF0} for the $R$-module $\Omega_{n}^F(M)$ , we get $\gd(\Omega_{n}^F(M))=0$; and hence $\gd(M)<\infty$. Indeed, by \cite[Corollary 2.4.8]{CF0}, $\gd(M)\leq\dim R$. Therefore, we could prove that every finitely generated $R$-modules $M$ satisfying the condition {$\sup\{~ i~|~\Ext_{R}^i(M,R)\neq0\}=\depth R-\depth M$} has finite $G$-dimension. But, this is a contradiction. Because, there exists a local Artinian ring $R$ which admits a finitely generated $R$-module $M$ satisfying the condition $\Ext_{R}^i(M,R)=0$ for all $i>0$, but $\gd(M)=\infty$ (see \cite{JSE}).}
\end{rem}

Considering the above mentioned remark, in the following proposition, we generalize the result \cite[Theorem 1.11]{KHY} which indicates that, over a Cohen-Macaulay local ring $R$ with a dualizing module $\Omega_{R}$, if $M$ is a finitely generated $R$-module of finite $G$-dimension, then $M$ is Cohen-Macaulay if and only if $M\otimes_{R}\Omega_{R}$ is Cohen-Macaulay.
\begin{prop}Let $(R,\fm)$ be a Cohen-Macaulay local ring with a dualizing $R$-module $\Omega_{R}$. Suppose that $M$ is a finitely generated $R$-module such that \emph{$\sup\{~ i~|~\Ext_{R}^i(M,R)\neq0\}=\depth R-\depth M$} and that \emph{$\Tor^R_{i}(M,\Omega_{R})=0$} for all $i>0$, then $M\otimes_{R}\Omega_{R}$ is Cohen-Macaulay if and only if $M$ is Cohen-Macauly.
\end{prop}
\begin{proof}First notice that, we may assume that $R$ is complete. Then $\Omega_{R}\cong\D_{\fm}(R)$. Also, since $R$ is Cohen-Macaulay, one gets $\gr(M)=\depth R-\dim M$. Therefore, the assertion follows from Theorem 5.6(ii).
\end{proof}
{$\mathbf{Acknowledgements}.$} I am very grateful to Professor Hossein Zakeri, Professor Kamran Divaani-Aazar and Professor Olgur Celikbas for their kind comments and assistance in the preparation of this manuscript.

\end{document}